\documentclass[12pt]{amsart}
\usepackage{verbatim}
\usepackage{amsmath}
\usepackage{enumerate}
\usepackage{amssymb}
\usepackage{color}
\usepackage{url}
\usepackage{graphicx}
\usepackage{fullpage}

\newcommand\+{\;\lower\plusheight\hbox{$+$}\;}
\newcommand\lldots{\;\lower\plusheight\hbox{$\cdots$}\;}

\newcommand{\Z}{\mathbb{Z}}
\newcommand{\Q}{\mathbb{Q}}

\newcommand{\C}{\mathbb{C}}

\newtheorem{Theorem}{Theorem}[section]
\newtheorem{Lemma}[Theorem]{Lemma}
\newtheorem{Corollary}[Theorem]{Corollary}

\newtheorem{Remark}[Theorem]{Remark}

\setbox0=\hbox{$+$}
\newdimen\plusheight
\plusheight=\ht 0 \setbox0=\hbox{$-$}
\newdimen\minusheight
\minusheight=\ht0

\setbox0=\hbox{$\cdots$}
\newdimen\cdotsheight
\cdotsheight=\plusheight

\begin{document}

\title{Sign-change of the Fourier coefficients of a Hauptmodul for $\Gamma_{0}(2)$}
\author{Bingyang Hu and Dongxi Ye$^\dag$}
\address{
Department of Mathematics, University of Wisconsin\\
480 Lincoln Drive, Madison, Wisconsin, 53706 USA.}
\email{bhu32@wisc.edu, lawrencefrommath@gmail.com}
\date{\today}
\maketitle
\begin{abstract}
In this short note, we aim to prove that the Fourier coefficients of the modular function
$
\frac{\eta^{24}(\tau)}{\eta^{24}(2\tau)}
$
possess a sign-change property.
\end{abstract}
\noindent
\numberwithin{equation}{section}
\allowdisplaybreaks


\section{Introduction}
\label{intro}


Let $\mathbb{H}$ be the upper half plane of $\C$, i.e., the set of complex numbers with positive imaginary part. Let $\Gamma$ denote the full modular group $\textrm{SL}_{2}(\mathbb{Z})$, i.e., the group of 2-by-2 matrices over integers of discriminant~1, and let $j(\tau)$ be the well known modular $j$-invariant. As one of the most famous functions in number theory, the modular $j$-invariant possesses numerous interesting properties. For example, it generates the function field of $\textrm{SL}_{2}\backslash\mathbb{H}^{*}$, where $\mathbb H^*=\mathbb H \cup \Q \cup \{\infty\}$, and its value at an imaginary quadratic point $\tau$, which is called a singular modulus, generate some ring class field of the imaginary quadratic field $\mathbb{Q}(\tau)$ \cite{Sh}. It is no exaggeration to say that the modular $j$-invariant appears everywhere in the area of number theory related to modular forms. Even in the recent active and developing study of sign changes of the Fourier coefficients of modular forms, one can be told by Asai et al \cite{AKN} that the signs of the Fourier coefficients of $\frac{1}{j(\tau)}$ are alternating. Those interesting properties of the modular $j$-invariant has motivated developments of a great variety of studies, of which the famous one is the traces of singular moduli which can be defined as follows.
Let $d$ be a positive integer congruent to $0$ or $3$ modulo $4$. We denote by $Q_d$ the set of positive-definite binary forms $Q(X, Y)=[a, b, c]=aX^2+bXY+cY^2, a, b, c \in \Z$ of discriminant $-d$, with usuall action of $\Gamma$. To each $Q \in Q_d$, we associate its unique root $\alpha_Q \in \mathbb H$ of $Q(X, 1)$ and $w_Q:=|\bar{\Gamma}_{Q}|$, where $\bar{\Gamma}_{Q}$ is the group of automorphisms of $Q$, i.e., the stabilizer of $Q$. Then the trace of singular moduli of discriminant $-d$ is defined by
$$
t(d)=\sum_{Q\in Q_{d}/\Gamma}\frac{j(\alpha_{Q})-744}{|\bar{\Gamma}_{Q}|}
$$
with the convention $t(0)=2, t(-1)=-1$ and $t(d)=0$ for $d<-1$ or $d \equiv 1, 2 \pmod4$. The traces of singular moduli play a key role in recovering some arithmetic information encoded in the Fourier coefficients of modular forms. For example, Kaneko \cite{K} showed that
$$
a(n)=\frac{1}{n}\left\{\sum_{r\in\mathbb{Z}}t(n-r^{2})+\sum_{\substack{r\geq1\\r\,odd}}\left((-1)^{n}t(4n-r^{2})-t(16n-r^{2})\right)\right\}
$$
for $n\geq1$, where $a(n)$ are the Fourier coefficients of the modular $j$-invariant. Inspired by Kaneko's work, Ohta \cite{O} derived analogous formulas for some genus zero modular functions such as the modular function
\begin{align*}
\frac{\eta^{24}(\tau)}{\eta^{24}(2\tau)}&=\frac{1}{q}\prod_{n=1}^{\infty}\frac{(1-q^{n})^{24}}{(1-q^{2n})^{24}}\\
&=\frac{1}{q}-24+276q-2048q^2+11202q^3-49152q^4+184024q^5-614400q^6+O(q^{7})\\
&=\sum_{n=-1}^{\infty}c(n)q^{n},
\end{align*}
where $\eta(\tau)$ is the Dedekind eta function defined by
$$
\eta(\tau)=q^{1/24}\prod_{n=1}^{\infty}(1-q^{n})
$$
with $q=\exp(2\pi i\tau)$, and he obtained that
\begin{align*}
c(n)
&=\frac{1}{n}\left\{\sum_{r\in\mathbb{Z}}t(n-r^{2})+\sum_{\substack{r\geq1\\r\,\,odd}}(-1)^{n}t(4n-r^{2})+24\sum_{\substack{d|n\\d\,\,odd}}d\right\}.
\end{align*}
This modular function has certain properties analogous to the modular $j$-invariant. For example, it is a generator for the function field of $\Gamma_{0}(2)\backslash\mathbb{H}^{*}$, where $\Gamma_0(2)$ denotes the congruence subgroup of level $2$, i.e.
$$
\Gamma_0(2)=\left\{
\begin{bmatrix}
    a  & b \\
    c  & d  \\
\end{bmatrix}
\in \Gamma\bigg| c \equiv 0  \pmod2 
\right\}.
$$
Also if one observes the Fourier coefficients $c(n)$ of $\frac{\eta^{24}(\tau)}{\eta^{24}(2\tau)}$, one may note that the signs of the first eight $c(n)$'s are alternating, like that of $\frac{1}{j(\tau)}$.
Recently, motivated by Kaneko's and Ohta's work, Matsusaka and Osanai \cite[Theorem 1.3]{MO} extend Ohta's formula for $c(n)$ to the case involving certain generalized traces of singular moduli, and moreover derive an amazing asymptotic formula for $c(n)$, namely, as $n\to\infty$,
\begin{equation} \label{asym}
c(n)\sim\frac{e^{2\pi\sqrt{n}}}{2n^{3/4}}\times
\begin{cases}
-1, & \mbox{if $n\equiv0\pmod{2}$};\\
1,  & \mbox{if $n\equiv1\pmod{2}$}
\end{cases}
\end{equation}
Clearly, as a consequence of the asymptotic formula \eqref{asym}, we note that the coefficient $c(n)$ possesses a sign-change property for large $n$. In this paper, 
we extend the above observations and get a full range description of the oscillatory behavior of $c(n)$. Namely, we have the following result.

\begin{Theorem}
\label{main}
For all integer $n\geq-1$, we have $(-1)^{n+1}c(n)>0$, i.e.,
$$
\begin{cases}
c(n)>0, &\mbox{if $n$ is odd};\\
c(n)<0, &\mbox{if $n$ is even}.
\end{cases}
$$
\end{Theorem}

An immediate consequence of Theorem \ref{main} is the following.
\begin{Corollary}
The Fourier coefficients $c(n)$ never vanish.
\end{Corollary}


\section{Proof of Theorem \ref{main}}


We need the following lemmas.

\begin{Lemma}[Robin] \label{robin}
For $n\geq3$, 
$$
\sum_{d|n}d<e^{\gamma}n\log\log{n}+\frac{n}{\log\log{n}},
$$
where $\gamma$ is Euler's constant.
\end{Lemma}
(see, e.g, \cite[Theorem 2]{R}).

\begin{Remark}
Under the Riemann Hypothesis, the above estimation can be sharpened as follows:
$$
\sum_{d|n}d<e^{\gamma}n\log\log{n}.
$$
\end{Remark}

\begin{Lemma}[Choi, Kim and Lim]\label{ckl}
The trace of singular moduli, $t(d)$, has the following sign-change property.
\begin{equation}
\label{tdsign}
\begin{cases}
t(d)>0, &\mbox{if $d\equiv0\pmod{4}$};\\
t(d)<0, &\mbox{if $d\equiv 3\pmod{4}$}.
\end{cases}
\end{equation}
\end{Lemma}
(see, e.g., \cite[Theorem 1]{CKL}).

\begin{Lemma}[Choi, Kim and Lim] \label{ckl2}
The trace of singular moduli, $t(d)$, has a lower bound and an upper bound as follows. 
\begin{enumerate}
\item[(i).]{If $d\equiv0\pmod{4}$, then
$$
\exp(\pi\sqrt{d})-\frac{1}{2}(2\pi d)^{\frac{3}{2}}\exp\left(\frac{\pi}{3}\sqrt{d}\right)\leq t(d)\leq
\exp(\pi\sqrt{d})+\frac{1}{2}(2\pi d)^{\frac{3}{2}}\exp\left(\frac{\pi}{3}\sqrt{d}\right).
$$
}
\item[(ii).] {If $d\equiv3\pmod{4}$, then
$$
-\exp(\pi\sqrt{d})-\frac{1}{2}(2\pi d)^{\frac{3}{2}}\exp\left(\frac{\pi}{3}\sqrt{d}\right)\leq t(d)\leq
\frac{1}{2}(2\pi d)^{\frac{3}{2}}\exp\left(\frac{\pi}{3}\sqrt{d}\right)-\exp(\pi\sqrt{d}).
$$
}
\end{enumerate}
\end{Lemma}
(see, e.g., \cite[(31)]{CKL}).

\begin{Lemma}[Zagier] \label{zagier}
For all positive integer $n$, we have
$$
\sum_{|r|<2\sqrt{n}}t(4n-r^{2})=
\begin{cases}
-4, &\mbox{if $n$ is a square}; \\
2,  &\mbox{if $4n+1$ is a square};\\
0,  &\mbox{otherwise}.
\end{cases}
$$
\end{Lemma}
(see, e.g., \cite[Theorem 2]{Z}).

\begin{Lemma}[Ohta]\label{ohta}
For all positive integer $n$, we have
\begin{align*}
c(n)
&=\frac{1}{n}\left\{\sum_{r\in\mathbb{Z}}t(n-r^{2})+\sum_{\substack{r\geq1\\r\,\,odd}}(-1)^{n}t(4n-r^{2})+24\sum_{\substack{d|n\\d\,\,odd}}d\right\}.
\end{align*}
\end{Lemma}
(see, e.g., \cite[Theorem 2.1]{O}).

\bigskip

\begin{proof}[\textbf{Proof of Theorem \ref{main}}]
Note that it suffices to prove the sign-change propety for the sequence $\{nc(n)\}_{n \ge -1}$, which by Lemma \ref{ohta}, is the same as
$$
\left\{\sum_{r\in\mathbb{Z}}t(n-r^{2})+\sum_{\substack{r\geq1\\r\,\,odd}}(-1)^{n}t(4n-r^{2})+24\sum_{\substack{d|n\\d\,\,odd}}d \right\}_{n \ge -1}.
$$


We divide our proof into four cases.

\medskip

\textit{Case I: $n \equiv 0 \pmod 4$.}

\medskip

We need to show that $nc(n)<0$ in this case. Suppose $n=4k$ for some $k \ge 0$, then we have
\begin{align*}
&\sum_{r\in\mathbb{Z}}t(4k-r^{2})+\sum_{\substack{r\geq1\\r\, odd}}t(16k-r^{2})+24\sum_{\substack{d|k\\2\nmid d}}d\\
&=\sum_{|r|\leq\sqrt{4k}}t(4k-r^{2})+\sum_{r\geq0}t(16k-r^{2})-\sum_{\substack{r\geq0\\r\,even}}t(16k-r^{2})+24\sum_{\substack{d|k\\2\nmid d}}d\\
& \quad \textrm{(By the definition of $t(d)$)} \\
&=\sum_{|r|\leq\sqrt{4k}}t(4k-r^{2})+\frac{1}{2}\sum_{|r|\leq\sqrt{16k}}t(16k-r^{2})+\frac{t(16k)}{2}-\sum_{\substack{r\geq0\\r\,even}}t(16k-r^{2})+24\sum_{\substack{d|k\\2\nmid d}}d\\
&=\sum_{|r|\leq\sqrt{4k}}t(4k-r^{2})+\frac{1}{2}\sum_{|r|\leq\sqrt{16k}}t(16k-r^{2})+\left(\frac{t(16k)}{2}-t(16k) \right)-\sum_{\substack{r>1 \\r\,even}}t(16k-r^{2})+24\sum_{\substack{d|k\\2\nmid d}}d\\
&=\sum_{|r|\leq\sqrt{4k}}t(4k-r^{2})+\frac{1}{2}\sum_{|r|\leq\sqrt{16k}}t(16k-r^{2})-\frac{t(16k)}{2}-\sum_{\substack{r> 1\\r\,even}}t(16k-r^{2})+24\sum_{\substack{d|k\\2\nmid d}}d\\
\end{align*}

By Lemma \ref{zagier}, we know that
$$
\sum_{|r|<\sqrt{4k}}t(4k-r^{2})+\frac{1}{2}\sum_{|r|<\sqrt{16k}}t(16k-r^{2})\leq3.
$$
Together with $t(0)=2$, this implies that
$$
\sum_{|r|\leq\sqrt{4k}}t(4k-r^{2})+\frac{1}{2}\sum_{|r|\leq\sqrt{16k}}t(16k-r^{2})\leq9.
$$

Moreover, it is clear that
$$
\sum_{\substack{d|k\\2\nmid d}}d\leq\sum_{d|k}d.
$$
Thus, we have
\begin{align*}
&\sum_{|r|\leq\sqrt{4k}}t(4k-r^{2})+\frac{1}{2}\sum_{|r|\leq\sqrt{16k}}t(16k-r^{2})-\frac{t(16k)}{2}-\sum_{\substack{r>1\\r\,even}}t(16k-r^{2})+24\sum_{\substack{d|k\\2\nmid d}}d\\
&\leq 9-\sum_{\substack{r>1\\r\,even}}t(16k-r^{2})-\frac{t(16k)}{2}+24\sum_{d|k}d.
\end{align*}
Since $16k-r^{2}\equiv0\pmod{4}$ when $r$ is even, then by Lemma \ref{ckl}, we know that
$$
\sum_{\substack{r>1\\r\,even}}t(16k-r^{2})>0.
$$

\textbf{Claim:}
When $k \ge 2$, we have
$$
9+24\sum_{d|k}d-\frac{t(16k)}{2}<0.
$$

Indeed, by Lemma \ref{robin} and Lemma \ref{ckl2}, we have
$$
\sum_{d|k}d<e^{\gamma}k\log\log{k}+\frac{k}{\log\log{k}}
$$
and
$$
t(16k) \ge \exp\left(\pi\sqrt{16k}\right)-\frac{1}{2}(32\pi k)^{\frac{3}{2}}\exp\left(\frac{\pi}{3}\sqrt{16k} \right).
$$
Hence, it suffices to show that when $k \ge 3$, the following inequality holds,
$$
9+24\left(e^{\gamma}k\log\log{(k)}+\frac{k}{\log\log{k}}\right)<\frac{1}{2}\left(\exp\left(\pi\sqrt{16k}\right)-\frac{1}{2}(32\pi k)^{\frac{3}{2}}\exp\left(\frac{\pi}{3}\sqrt{16k}\right)\right).
$$
However, the proof for the above inequality is elementary and hence we omit it here. Thus, we already show that for $k \ge 3$, $c(n)=c(4k)<0$.  Combining this estimation with the fact that $c(0)=-24$, $c(4)=-49152$ and $c(8)=-5373952$, we get the desired result. 

\medskip

\textit{Case II: $n \equiv 2 \pmod 4$.}

\medskip

Suppose $n=4k+2$ for some $k \ge 0$. Similarly, we start with Lemma \ref{ohta} to find
\begin{align*}
nc(n)&=\sum_{r\in\mathbb{Z}}t(4k+2-r^{2})+\sum_{\substack{r\geq1\\r \, odd}}t(16k+8-r^{2})+24\sum_{\substack{d|(2k+1)}}d.
\end{align*}
Note that $4k+2-r^{2}\equiv 1 \ \textrm{or} \ 2\pmod{4}$, and hence by the definition of $t(d)$, it follows that $t(4k+2-r^2)=0$ for all $r \in \Z$, which implies 
$$
\sum_{r \in \Z} t(4k+2-r^2)=0.
$$
By replacing the role of $n=4k$ in Case I by $n=4k+2$, we can deduce that
$$
\sum_{\substack{r\geq1\\r,\, odd}}t(16k+8-r^{2})+24\sum_{\substack{d|(2k+1)}}d<0
$$
for all $k \ge 2$. Combining this estimation with the fact that $C(2)=-2048$, we get $c(n)=c(4k+2)<0$.

Thus, from the previous two cases, we see that $c(n)<0$ for $n$ even.

\medskip

\textit{Case III: $n \equiv 1 \pmod 4$.}

\medskip
 
Let $n=4k+1$ for some $k \ge 1$. In this case, we have
$$
nc(n)=\sum_{r\in\mathbb{Z}}t(4k+1-r^{2})-\sum_{\substack{r\geq1\\r\, odd}}t(16k+4-r^{2})+24\sum_{\substack{d|(4k+1)}}d.
$$
Since $4k+1-r^{2}\equiv 0 \ \textrm{or} \ 1\pmod{4}$, then by Lemma \ref{ckl} together with $t(d)=0$ if $d\equiv1,2\pmod{4}$, we have
$$
\sum_{r\in\mathbb{Z}}t(4k+1-r^{2})>0.
$$
Moreover, for $r$ odd, we have $16k+4-r^{2}\equiv3\pmod{4}$. Using Lemma \ref{ckl} again, we find that
$$
\sum_{\substack{r\geq1\\r\, odd}}t(16k+4-r^{2})<0.
$$
 Therefore, $nc(n)>0$ for $n=4k+1$.

\medskip

\textit{Case IV: $n \equiv 3 \pmod 4$.}

\medskip
 
Finally, for the case $n=4k+3, k \ge -1$, again, we start with
\begin{align*}
nc(n)&=\sum_{r\in\mathbb{Z}}t(4k+3-r^{2})-\sum_{\substack{r\geq1\\r \, odd}}t(16k+12-r^{2})+24\sum_{\substack{d|(4k+3)}}d\\
&=\sum_{|r|\leq\sqrt{4k+3}}t(4k+3-r^{2})-\sum_{\substack{1\leq r\leq\sqrt{16k+12}\\r \, odd}}t(16k+12-r^{2})+24\sum_{\substack{d|(4k+3)}}d\\
&= \left( \sum_{|r|\leq\sqrt{4k+3}}t(4k+3-r^{2})-t(16k+11) \right)-\sum_{\substack{2\leq r\leq\sqrt{16k+12}\\r \, odd}}t(16k+12-r^{2})+24\sum_{\substack{d|(4k+3)}}d\\
\end{align*}
Clearly, if $r$ is odd, $16k+12-r^{2}\equiv3\pmod{4}$, then $t(16k+12-r^{2})<0$, and hence
$$
\sum_{\substack{2\leq r\leq\sqrt{16k+12}\\r \, odd}}t(16k+12-r^{2})<0.
$$
While for $r \in \Z$, $4k+3-r^{2}\equiv 2 \ \textrm{or} \ 3 \pmod{4}$, then by Lemma \ref{ckl}, we have
$$
 \sum_{|r|\leq\sqrt{4k+3}}t(4k+3-r^{2})=\sum_{\substack{|r|\leq\sqrt{4k+3}\\r \, even}}t(4k+3-r^{2}).
$$
Thus, by Lemma \ref{ckl2}, (ii) and an easy calculation, we see that for $k \ge 3$, 
\begin{eqnarray*}
&& \sum_{|r|\leq\sqrt{4k+3}}t(4k+3-r^{2})-t(16k+11)\\
&& \ge \exp \left(\pi \sqrt{16k+11} \right)-\frac{\left[2\pi(16k+11)\right]^{\frac{3}{2}}}{2} \exp \left(\frac{\pi \sqrt{16k+11}}{3} \right)\\
&& \quad -2 \sqrt{4k+3} \left( \exp\left(\pi \sqrt{4k+3}\right)+\frac{\left[2\pi(4k+3)\right]^{\frac{3}{2}}}{2} \exp \left( \frac{\pi \sqrt{4k+3}}{3} \right) \right)\\
&&>0.
\end{eqnarray*}
Thus, we have seen that when $n=4k+3, k \ge 3$, $c(n)>0$. Combining this estimation with the fact that $c(-1)=1, c(3)=11202, c(7)=1881471$ and $c(11)=91231550$, it follows that $c(n)>0$ for $n=4k+3$ with $k \ge 1$. 

Clearly, the arguments in Case III and Case IV imply that $c(n)>0$ for $n$ odd.

\begin{Remark}
For general cases, say $p\geq3$ and $(p-1)|24$, we will need to consider the generalized traces of singular moduli (see, e.g., \cite{CJKK}), and these will be treated in  subsequent work of the authors \cite{HY}.
\end{Remark}


\end{proof}


\begin{thebibliography}{99}

\bibitem{AKN}{T. Asai, M. Kaneko, and H. Ninomiya, \emph{Zeros of certain modular functions and an application}, Comment. Math. Univ. St. Paul., {\bf46} (1997), 93--101.}

\bibitem{CJKK}{D. Choi, D. Jeon, S.-Y. Kang, and C. H. Kim, \emph{ Traces of singular moduli of arbitrary level modular functions}, Int. Math. Res. Not. IMRN 2007, no. 22, Art. ID rnm110, 17 pp.}

\bibitem{CKL}{D. Choi, B. Kim and S. Lim, \emph{Sign-Periodicity of Traces of Singular Moduli}, Mathematics, {\bf1} (2013), 111--118}

\bibitem{HY}{B. Hu and D. Ye, \emph{Sign changes of Fourier coefficients of Hauptmoduls for genus zero groups}, in preparation.}

\bibitem{K}{M. Kaneko, \emph{The Fourier coefficients and the singular moduli of the elliptic modular
function $j(\tau)$}, Mem. Fac. Engrg. Design Kyoto Inst. Tech. Ser. Sci. Tech. {\bf44} (1995), 1--5.}

\bibitem{MO}{T. Matsusaka and R. Osanai, \emph{Arithmetic formulas for the Fourier coefficients of Hauptmoduln of level 2, 3, and 5}, Proc. Amer. Math. Soc., {\bf145} (2017), 1383--1392.}

\bibitem{O}{K. Ohta, \emph{Formulas for the Fourier coefficients of some genus zero modular functions}, Kyushu J. Math. {\bf63} (2009), 1--15.}

\bibitem{R}{G. Robin, \emph{Grandes valeurs de la fonction somme des diviseurs et hypoth$\mbox{$\grave{e}$}$se de Riemann},  J. Math. Pures Appl. {\bf63} (1984), 187--213.}

\bibitem{Sh}{G. Shimura, \emph{Introduction to the Theory of Automorphic Forms}, Iwanami Publishing Company and Princeton University Press, 1974.}


\bibitem{Z}{D. Zagier, \emph{Traces of singular moduli}, In \emph{Motives, Polylogarithms and Hodge Theory}; Int. Press Lect. Ser. 3; International Press: Somerville, MA, USA, 2002; 209--244.}
\end{thebibliography}
\end{document}